\documentclass[11pt]{amsart}
\usepackage{amsmath,amssymb,epsfig,color, amsthm, mathrsfs}
\usepackage{wrapfig,lipsum,floatrow,sidecap}
\usepackage{graphicx}
\usepackage{multirow}
\usepackage{hhline}
\usepackage{url}
\graphicspath{ {images/} }
\allowdisplaybreaks

\newtheorem{theorem}{Theorem}[section]

\newtheorem{corollary}[theorem]{Corollary}
\newtheorem{definition}[theorem]{Definition}
\newtheorem{lemma}[theorem]{Lemma}
\newtheorem{proposition}[theorem]{Proposition}
\newtheorem{remark}[theorem]{Remark}
\newtheorem{example}[theorem]{Example}

\newcommand{\vanish}[1]{}\parskip=12pt

\def\p{\prime}
\def\pp{{\prime\prime}}
\def\x{\textbf{x}}
\def\b{\textbf{b}}

\def\G{\overline{G}}
\def\R{\mathbf{R}}

\def\A{\mathcal{A}}
\def\B{\backslash}
\def\C{\mathcal{C}}
\def\DD{\mathcal{D}}
\def\D{\mathbb{D}}
\def\B{\mathcal{B}}

\def\L{\mathcal{L}}
\def\K{\mathcal{K}}

\numberwithin{equation}{section}

\begin{document}
\title{Anti-parallel links as boundaries of knotted ribbons}
\author{Yuanan Diao}
\address{Department of Mathematics and Statistics\\
University of North Carolina Charlotte\\
Charlotte, NC 28223}
\email{ydiao@uncc.edu}

\begin{abstract}
We consider a two component link $\mathcal{L}$ which is the boundary of a knotted ribbon with knot type $\mathcal{K}$. If the two components of $\mathcal{L}$ are assigned opposite (anti-parallel) orientations and $\mathcal{K}$ is a special alternating knot, namely a knot with a reduced alternating knot diagram in which the crossings are all positive or all negative, we show that the braid index of $\mathcal{L}$ is bounded below by ${\rm{Cr}}(\mathcal{K})+2$, where ${\rm{Cr}}(\mathcal{K})$ is the minimum crossing number of $\mathcal{K}$. A long standing open conjecture states that the ropelength of any alternating knot is bounded below by its minimum crossing number multiplied by a positive constant. Using the above result, we are able to prove that this conjecture holds for the special alternating knots.
\end{abstract}

\maketitle
\section{Introduction}\label{s1}
A knotted ribbon with a given knot type $\K$ is an embedding of the annulus $\{1\le x^2+y^2\le 2:\ x, y\in \R\}$ into $\R^3$ such that the centerline (axis) of the annulus is a knot of the given knot type $\K$. The boundary of a knotted ribbon with knot type $\K$ is a two component link denoted by either $L_p(\K)$ or $L_a(\K)$, depending on whether the two components are assigned parallel or anti-parallel orientations respectively. 
In \cite{White} J.~White proved the famous formula 
\begin{equation}\label{linking_eq}
{\rm{Lk}}(L_p(\K))={\rm{Tw}}((L_p(\K))+{\rm{Wr}}((L_p(\K)),
\end{equation} 
where ${\rm{Lk}}(L_p(\K))$ is the linking number of $L_p(\K)$ (a link invariant), ${\rm{Tw}}((L_p(\K))$ is the twist of $L_p(\K)$ and 
${\rm{Wr}}((L_p(\K))$ is the writhe of the centerline (or the axis) of $L_p(\K)$.  Equation (\ref{linking_eq}) has a very important application in DNA research since a knotted ribbon provided an ideal mathematics model for circular double stranded DNA in which the two sugar-phosphate  backbones of the double helix are modeled by the boundary curves of the knotted ribbon \cite{White2}. Inspired by this work, in this paper we study another link invariant related to a link that is the boundary of a knotted ribbon, namely the braid index. The braid index of $L_p(\K)$ can be easily shown to be $2\b(\K)$ where $\b(\K)$ is the braid index of $\K$, thus the case of interest to us is the braid index of $L_a(\K)$. 

\begin{figure}[!hbt]
\begin{center}
\includegraphics[scale=0.8]{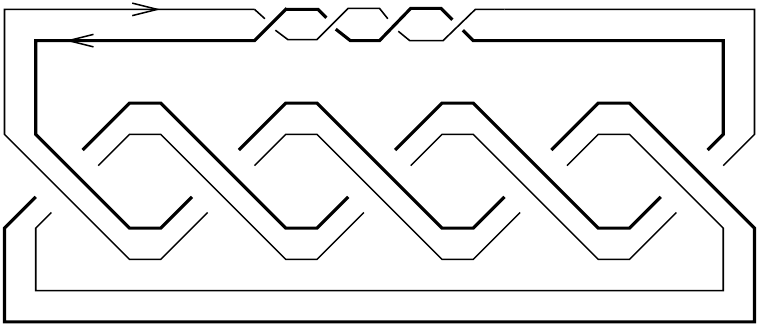}
\end{center}
\caption{The projection diagram of an embedded annulus whose two components are equivalent to the $5_1=T(2,5)$ torus knot. The top portion of the diagram contains 2 full positive twists and the linking number of the link is thus $(-10+4)/2=-3$.\label{double}}
\end{figure}

\medskip
The main result of this paper is that if $\K$ is an alternating knot with a reduced knot projection in which the crossings are either all positive or all negative (such an alternating knot is called a {\em special alternating knot}  by Murasugi in \cite{Mu1960}), then the braid index of $L_a(\K)$ is 
bounded below both by the minimum crossing number of $\K$ and the absolute linking number of $L_a(\K)$ (multiplied by some positive constant). More specifically, let $\K$ be an alternating knot with a reduced knot projection $D$ in which the crossings are either all positive or all negative, then $L_a(\K)$ admits a projection diagram $\D_k$ as depicted in Figure \ref{double}, where the two components of $L_a(\K)$ are realized as $D$ and a parallel copy of $D$ (but with opposite orientations). A crossing in $D$ becomes a 4 crossing junction in $\D_k$. The additional $|k|$ full twists between these two components are determined by $k=w(D)+{\rm{Lk}}(L_a(\K))$, where $w(D)$ is the writhe of the diagram $D$ (namely the sum of the number of positive crossings minus the number of negative crossings in $D$). $k>0$ and $k<0$ indicate that the crossings in these full twists are positive or negative respectively. Let $c(D)$ be the number crossings in $D$ and $s(D)$ be the number of Seifert circles (to be discussed in the next section), then our main result can be formulated as the following theorem.

\begin{theorem}\label{main_theorem}
Let $\K$ be a special alternating knot with a reduced alternating link diagram such that $c(D)=|w(D)|$, then the braid index $\b(\D_k)$ of $\D_k$, which is the same as that of $L_a(\K)$ with $k=w(D)+{\rm{Lk}}(L_a(\K))$, is bounded below by
$
c(D)+2+\rho_k(D)
$
where
$$
\rho_k(D)=\left\{
\begin{array}{ll}
\max\{0,|k|+s(D)-c(D)-2\}& \ {\rm if}\ w(D)>0, \ k\le 0\ {\rm or}\ w(D)<0, \ k\ge 0;\\
\max\{0,|k|-s(D)\}& \ {\rm if}\ w(D)>0, \ k\ge 0\ {\rm or}\ w(D)<0, \ k\le 0.
\end{array}
\right.
$$
\end{theorem}

Since $\rho_k(D)$ grows linearly in terms of $|{\rm{Lk}}(L_a(\K))|$ (after it passes some threshold), the braid index of $L_a(\K)$ is bounded below by $c(D)+2$ plus a linear term of $|{\rm{Lk}}(L_a(\K))|$. This result has a significant implication concerning the ropelength problem of knots as explained below.

\medskip
The ropelength of a link is  defined (intuitively) as the minimum length of unit thickness ropes needed to tie the link. Let $\K$ be an un-oriented link, $Cr(\K)$ be the minimum crossing number of $\K$ and $R(\K)$ be the ropelength of $\K$. A fundamental question in geometric knot theory asks how $R(\K)$ is related to $Cr(\K)$. In general the determination of the precise ropelength of a non-trivial link is a difficult problem and no known precise formula exists for $R(\K)$ even for the simplest nontrivial knot, namely the trefoil. However we have gained much knowledge about the asymptotic behavior of $R(\K)$ in terms of $Cr(\K)$. For example,  it is known that $a_0(Cr(\K))^{3/4}\le R(\K)\le a_1 Cr(\K)\ln^5(Cr(\K))$ for some positive constants $a_0$ and $a_1$ for any $\K$ \cite{Buck, Buck2,Diao2019_2}, and for any power  $p$ that is between $3/4$ and $1$, there exist families of infinitely many links such that the ropelength of links from these families grows proportionally to $(Cr(\K))^{p}$ \cite{Cantarella1998,Diao2020,Diao1998,Diao2003}. It has been conjectured that the ropelength of any alternating link is at least proportional to its crossing number. While this conjecture has been shown to be true for alternating links whose maximum braid indices are  proportional to their crossing numbers, it is still wide open in general. Using the main result of this paper, we are able to prove that this conjecture holds for all special alternating knots regardless whether their braid indices are proportional to their crossing numbers. 

\medskip
We arrange the rest of the paper as follows. In the next section, we introduce the Seifert circle decomposition of oriented link diagrams and their corresponding Seifert graphs. Some graph theory terminology specific to this paper will also be introduced. In Section \ref{MP_sec} we introduce the Murasugi-Przytycki reduction operation. Some important results regarding the HOMFLY-PT polynomial are given in Section \ref{Homfly_sec} as the HOMFLY-PT polynomial is our main tool in estimating the braid index. Section \ref{proof_sec} is devoted to the proof of Theorem \ref{main_theorem} and contains the bulk of the technical details of this paper. In the last section we apply \ref{main_theorem} to the ropelength problem of knots.

\section{Seifert circle decomposition and Seifert graphs}

\subsection{Seifert circle decomposition} 
\begin{figure}[htb!]
\includegraphics[scale=.9]{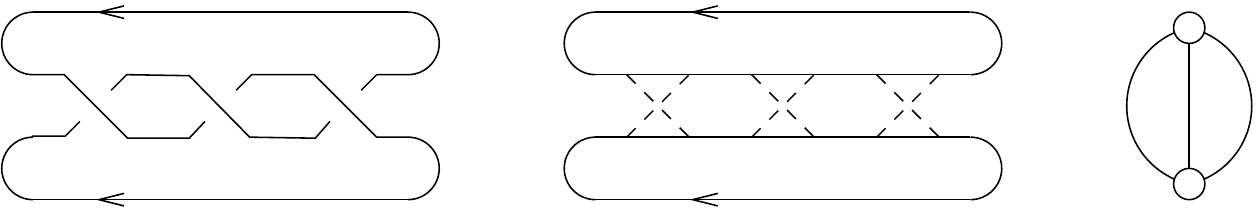}
\caption{Left: $D$ is the diagram of a positive trefoil; Middle: The Seifert circle decomposition of $D$; Right: The Seifert graph $G_{S}(D)$.}
\label{trefoil}
\end{figure}

The Seifert circle decomposition of $D$ is the collection of disjoint topological circles (called Seifert circles or $s$-circles for short) obtained after all crossings in $D$ are smoothed (as shown in Figure \ref{trefoil}). Notice that since $D$ is positive and is alternating, it can be drawn in a way that no $s$-circle can contain other $s$-circles in its interior.

\subsection{Trapped $s$-circles} 
\begin{figure}[htb!]
\includegraphics[scale=.9]{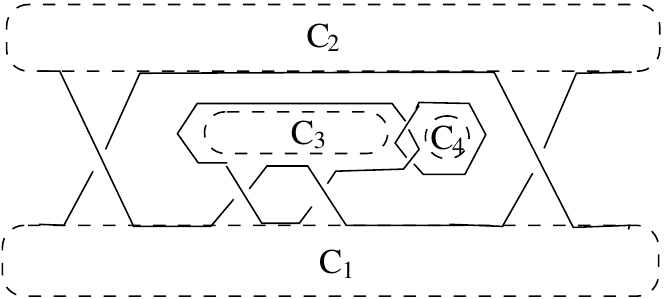}
\caption{$C_3$ and $C_4$ are trapped by $C_1$ and $C_2$. }
\label{trap}
\end{figure}

Consider three $s$-circles $C_1$, $C_2$ and $C_3$ as shown in Figure \ref{trap} where $C_3$ is bounded within the
 topological disk created by arcs of $C_1$, $C_2$ and the two consecutive
    crossings as shown in the figure. We say that $C_3$ is {\em trapped} between $C_1$ and $C_2$. Similarly, $C_4$ is  also trapped between $C_1$ and $C_2$.  Since $D$ is a special knot diagram in our case, the interior of any $s$-circle does not contain any other $s$-circle. Thus we can always use flype moves to free any trapped $s$-circles, hence we can assume that $D$ is free of trapped $s$-circles from this point on. An immediate consequence of this is that if two $s$-circles $C_1$ and $C_2$ share crossings, then as we travel along either $C_{1}$ or $C_{2}$, we encounter all these crossings consecutively without running into crossings between them and other $s$-circles.
    
\subsection{Seifert graphs}\label{sgraphsection}
We assume that our reader has the basic knowledge of graph theory and will only define the concepts that are specific to this paper. The common graph theory terminology used in this paper can be found in standard textbooks such as \cite{West}. 

\begin{definition}{\em
Let $D$ be a reduced alternating diagram, then shrinking each $s$-circle to a vertex and changing each crossing between two $s$-circles to an edge connecting the two corresponding vertices, we obtain a (bipartite) plane graph. We denote this graph by $G_S(D)$ and call it the {\em Seifert graph} of $D$. See the right side of Figure \ref{trefoil} for an example.
}
\end{definition}

If two $s$-circles in $S(D)$ share one and only one crossing between them, then this crossing is called a {\em lone crossing}. Two $s$-circles in $S(D)$ sharing a lone crossing corresponds to two vertices in $G_S(D)$ connected by one and only one edge which we shall call a {\em lone edge}. Notice that a lone edge in $G_S(D)$ cannot be a bridge, since otherwise it would correspond to a nugatory crossing of $D$, but that is not possible since $D$ is reduced. 
Sometimes it is more convenient to use a different version of the Seifert graph of $D$, denoted by $\G_S(D)$. $\G_S(D)$ is similar to $G_{S}(D)$ except that it is a simple graph with its edges weighted. An edge of weight $k$ in $\G_S(D)$ connecting two vertices $v_{1}$ and $v_{2}$ corresponds to the case when the same vertices $v_{1}$ and $v_{2}$ in $G_{S}(D)$ are connected by $k$ edges. A lone edge in $G_{S}(D)$ is an edge of weight one in  $\G_{S}(D)$. Since $K$ is a knot, $\G_{S}(D)$ cannot contain a bridge edge of even weight. 

\begin{figure}[htb!]
\includegraphics[scale=.5]{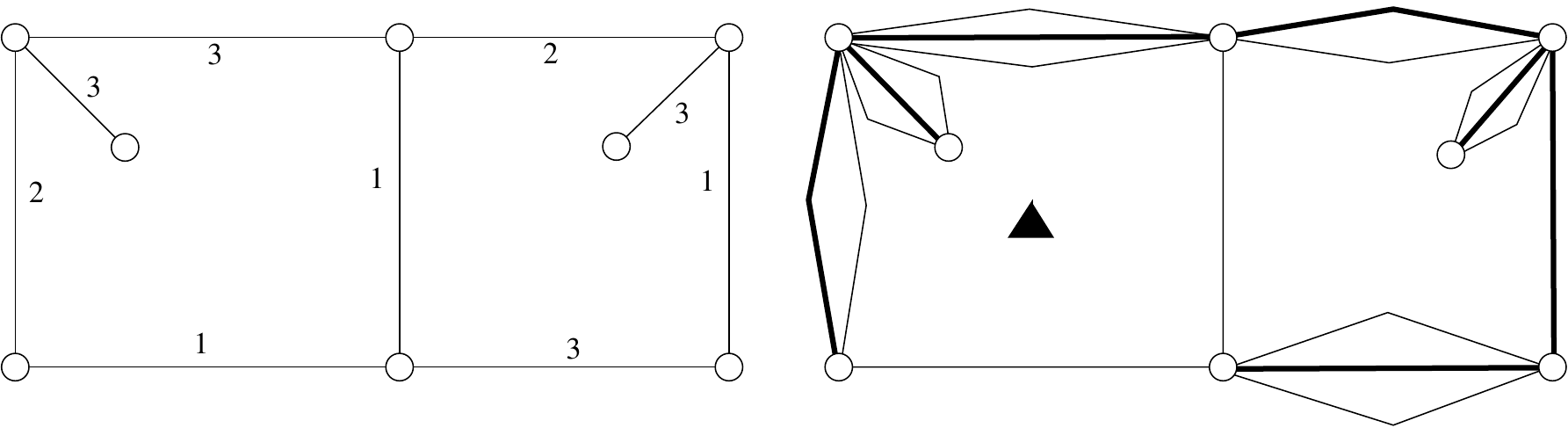}
\caption{Left: The edge weighted simple Seifert graph $\G_{S}(D)$ of a proper knot diagram $D$; Right: A good spanning tree of $G_{S}(D)$ with respect to the face marked by the triangle.}
\label{APknot}
\end{figure}

\begin{definition}{\em
A face of a plane graph $G$ is said to be {\em non-separating} if deleting the edges on its boundary does not disconnect the graph $G$. The graph $G$ is said to be {\em proper} if every face of it is non-separating. $D$ is said to be {\em proper} if every face in $G_{S}(D)$ is non-separating. See Figure \ref{APknot} for an example.
}
\end{definition}

\medskip
\begin{definition}{\em
A set $\Phi$ of edges in $\G_{S}(D)$ is said to be a {\em maximal 2-cut set} if it satisfies the following conditions: (i) $\Phi$ contains only lone edges and $|\Phi|\ge 2$; (ii) deleting any two edges in $\Phi$ disconnects $\G_{S}(D)$ and (iii) no other lone edges can be added to $\Phi$ in order for (ii) to hold.
}
\end{definition}

\begin{figure}[htb!]
\includegraphics[scale=.8]{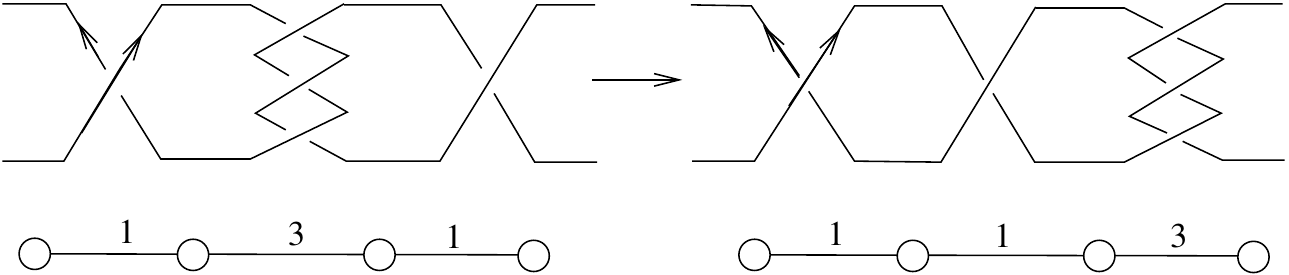}
\caption{Left: Part of the diagram $D$ with two lone crossings corresponding to two edges in a maximal 2-cut set of $\G_{S}(D)$. The shaded regions are isolated after the two lone crossings/edges are smoothed/deleted in $D$/$\G_{S}(D)$; Right: After  a  proper flype move, the two lone edges connect to a common vertex that is only connected by these two lone edges.}
\label{flype2}
\end{figure}

\begin{remark}\label{maximalpathremark}{\em
If $\G_{S}(D)$ contains a maximal 2-cut set $\Phi$, then we can perform flypes on $D$ so that in the resulting knot diagram $D^\p$, the lone crossings corresponding to the lone edges in $\Phi$ occur in a consecutive manner. In other word, if $D^\p$ is the new diagram, then the lone edges in $\Phi$ form a path in $\G_{S}(D^\p)$ such that the internal vertices of this path are all of degree 2.  Thus from this point on we will assume that a maximal 2-cut set of $\G_{S}(D)$, if it exists, consists of a path of lone edges whose internal vertices are all of degree 2. See Figure \ref{flype2} for an illustration of a simple case.
}
\end{remark}

\begin{lemma}
Let $D$ be a special alternating knot diagram, then the following conditions are equivalent.\\
\noindent
(i) $D$ is proper;\\
\noindent
(ii) For any face $F$ of $G_{S}(D)$, there exists a spanning tree $T_{F}$ of $G_{S}(D)$ such that no edges on $T_{F}$ are on the boundary of $F$.\\
\noindent
(iii) $\G_{S}(D)$ does not contain any maximal 2-cut set.
\end{lemma}

\begin{proof}
(i) $\Longrightarrow$ (ii) If $D$ is proper and $F$ is any face of $G_{S}(D)$, then deleting the edges on the boundary $\partial F$ of $F$ does not disconnect the graph,  which  means we can construct a spanning tree of $G_{S}(D)$ without using any edges on $\partial F$. (ii) $\Longrightarrow$ (iii) If $\G_{S}(D)$ contains a maximal 2-cut set $\Phi$, then it is necessary that the lone edges in $\Phi$ are all on the boundary of a face $F$. Since deleting any two lone edges in $\Phi$ will disconnect $\G_{S}(D)$, this means that if we delete the edges on the boundary of $F$, we will disconnect $G_{S}(D)$ as well. Thus it is not possible for us to construct a spanning tree of $G_{S}(D)$ without using edges on the boundary of $F$, which contradicts the given condition. (iii) $\Longrightarrow$ (i) Consider any face $F$ of $G_{S}(D)$. First consider the case that $\partial F$ contains a cycle $\Gamma$ in $G_{S}(D)$ which has length 2 and corresponds to an edge $\overline{\gamma}$ of weight at least 2 in $\G_{S}(D)$.  If $\overline{\gamma}$ is not a bridge edge of $\G_{S}(D)$, then deleting the edges in $\gamma$ obviously will not disconnect $G_{S}(D)$. If $\overline{\gamma}$ is  a bridge edge of $\G_{S}(D)$, then it is necessary that its weight is an odd integer that is at least 3 since $D$ is a reduced knot diagram. Therefore deleting the edges in $\gamma$ will not disconnect $G_{S}(D)$ either. 
If $\partial F$ contains a cycle $\gamma$ that is of length at least 4, then let us remove the edges of $\gamma$ one at a time. If the edge removed is not a lone edge, then obviously this does not disconnect $G_{S}(D)$. If the edge is the first lone edge encountered, removing the edge will not disconnect $G_{S}(D)$ either since the lone edge cannot be a bridge edge. Thus the only way $G_{S}(D)$ becomes disconnected is when we encounter a second lone edge in $\gamma$. But this means that a maximal 2-cut set exists which contradicts the given condition. Since the boundary of $F$ consists of cycles, this proves that removing the boundary of $F$ will not disconnect $G_{S}(D)$. 
\end{proof}

For any plane graph $G$, its number of faces is given by $e(G)-v(G)+2$, where $e(G)$ and $v(G)$ are the number of edges and vertices of $G$ respectively. In the case that $G=G_{S}(D)$, we have $e(G)=c(D)$ and $v(G)=s(D)$, where $c(D)$ and $s(D)$ are the number of crossings and number of $s$-circles in  $D$ respectively.

\section{The Murasugi-Przytycki reduction operation}\label{MP_sec}

\begin{figure}[htb!]
\includegraphics[scale=.7]{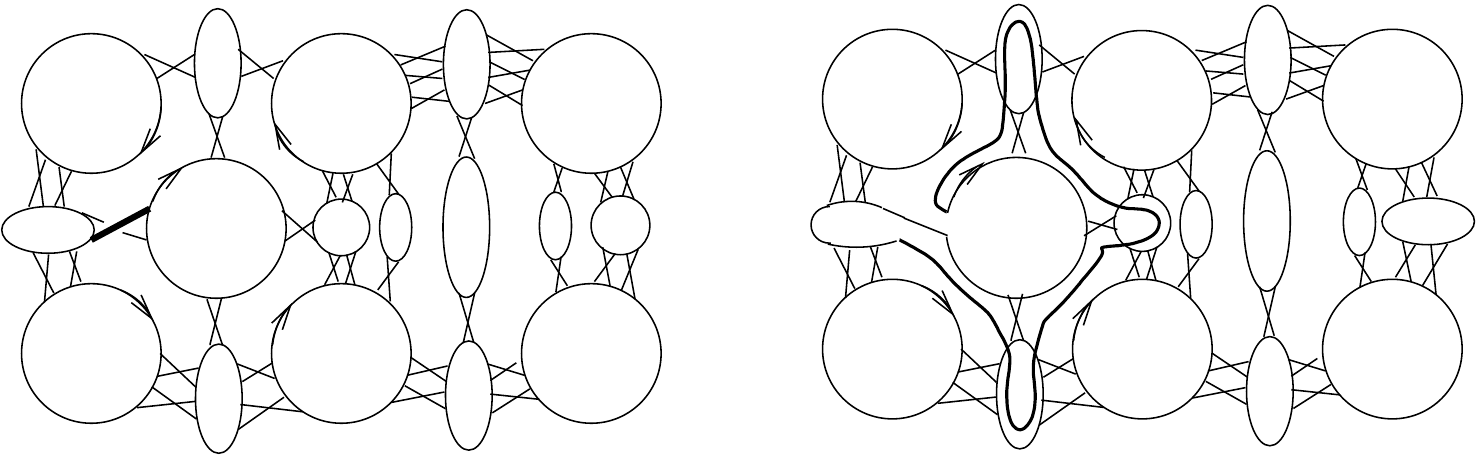}\\
\includegraphics[scale=.8]{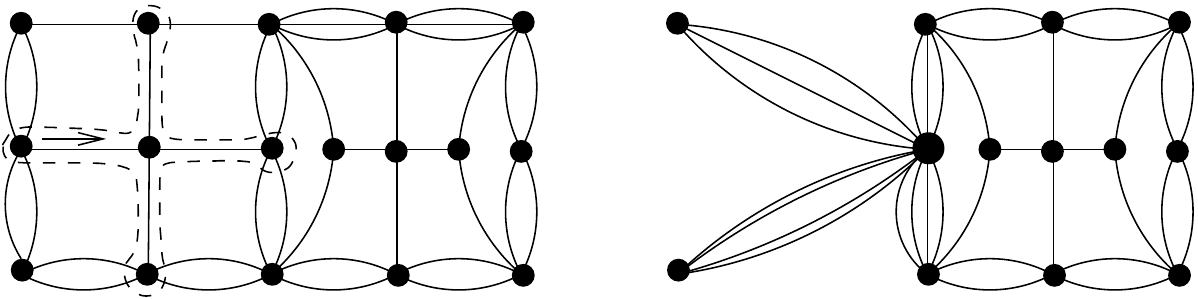}
\caption{How an overpass at a lone crossing is rerouted and the effect of the rerouting on the corresponding Seifert graph of the diagram. }
\label{reduction}
\end{figure}

Consider a link diagram $D$ with the property that none of its Seifert circles contains other Seifert circles in its interior. Here $D$ is more general, not necessarily a positive or negative alternating knot diagram. If $D$ has a lone crossing $x$ between two $s$-circles $C_1$ and $C_2$, then we can reroute the overpass at this crossing as shown in Figure \ref{reduction}.

The overpass to be rerouted is marked by a thick line in the top left diagram of Figure \ref{reduction} and the rerouted pass is marked by the thick line in the top right diagram (the thick line is the over strand at each crossing). The effect of this rerouted strand to the Seifert circle decomposition is that $C_1$ and $C_2$ are combined into one $s$-circle $C_0$, and any $s$-circle $C^\p$ sharing crossings with $C_2$ becomes an $s$-circle contained within this new $s$-circle $C_0$ (we say that $C^{\p}$ is {\em swallowed} by $C_{0}$), while any $s$-circle sharing crossings with $C^\p$ will now share the same crossings with $C_0$ instead. If we ignore the $s$-circles swallowed by $C_{0}$, then the effect of this operation on $G_{S}(D)$ is that the vertices $v_{1}$, $v_{2}$ corresponding to $C_1$, $C_2$, and any vertex $v$ corresponding to an $s$-circle sharing crossings with $C_{2}$, is contracted to the same  vertex $v_{0}$. The bottom of Figure \ref{reduction} shows this effect on the Seifert graph of the diagram. Notice that in the above described operation, we can also use $C_{1}$ in the place of $C_{2}$, which will lead to a different diagram and a different Seifert graph as it will contract the vertices connected with $v_{1}$ instead. This operation is known as the Murasugi-Przytycki reduction operation (we shall call it the MP operation for short) \cite{MP}.

\begin{definition}{\em
For a given link diagram $D$, we use $r^+(D)$ and $r^-(D)$ to denote the maximum number of MP operations that can be performed on positive and negative lone crossings respectively.
}
\end{definition}

\begin{remark}\label{reduction_number}{\em
If a lone crossing is positive (negative), then applying the MP operation to the diagram $D$ results in a new diagram $D^\p$ with $s(D^\p)=s(D)-1$ and $w(D^\p)=w(D)-1$ ($w(D^\p)=w(D)+1$), since the overpass crosses an existing Seifert circle an even number of times in the re-routing process hence does not change the writhe of $D$ except that the lone positive (negative) crossing is eliminated. Thus in general, $D$ can be deformed via the MP operations to new diagrams $D^\p$ and $D^\pp$ such that $s(D^\p)=s(D)-r^+(D)$, $w(D^\p)=w(D)-r^+(D)$ and $s(D^\pp)=s(D)-r^-(D)$, $w(D^\pp)=w(D)+r^-(D)$.
}
\end{remark}

\section{The HOMFLY-PT polynomial}\label{Homfly_sec}

\begin{figure}[htb!]
\includegraphics[scale=.5]{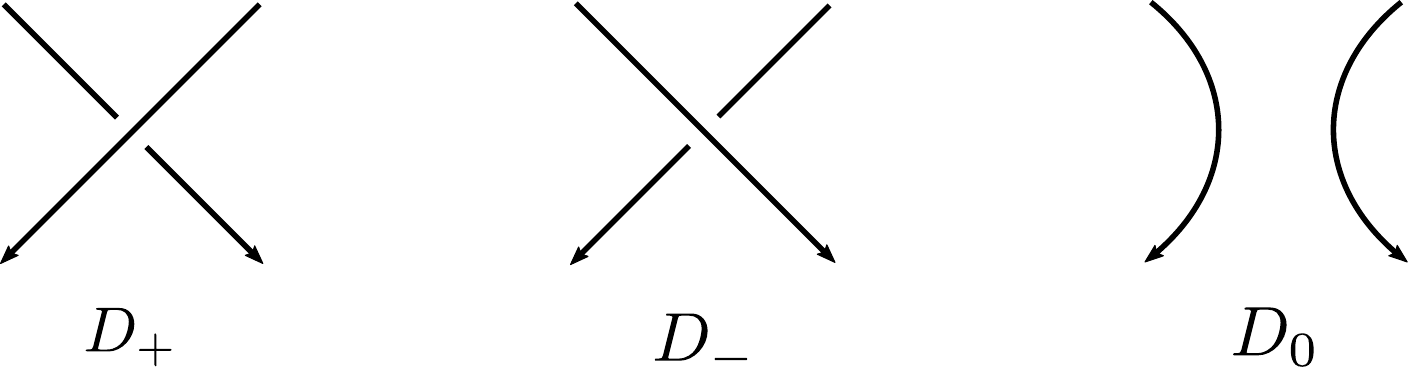}
\caption{The sign convention at a crossing of an oriented link and the smoothing of the crossing: the crossing in $D_+$ ($D_-$) is positive (negative) and is assigned $+1$ ($-1$) in the calculation of the writhe of the link diagram.}
\label{fig:cross}
\end{figure}

Let us recall that the HOMFLY-PT polynomial  of an oriented link is defined using any diagram $D$ of the link with the following rules: (a) If $D_1$ and $D_2$ represent the same link, then $H(D_1,z,a)=H(D_2,z,a)$; (b) $H(D,z,a)=1$ if $D$ is an unknot and (c) $aH(D_+,z,a) - a^{-1}H(D_-,z,a) = zH(D_0,z,a)$ where $D_+$, $D_-$ and $D_0$ are identical link diagrams except at one crossing as shown in Figure \ref{fig:cross}. For the sake of simplicity we shall use $H(D)$ for $H(D,z,a)$.

For our purposes, we will actually be using the following two equivalent forms of the skein relation (c) above:
\begin{eqnarray}
H(D_+,z,a)&=&a^{-2}H(D_-,z,a)+a^{-1}zH(D_0,z,a),\label{Skein1}\\
H(D_-,z,a)&=&a^2 H(D_+,z,a)-azH(D_0,z,a).\label{Skein2}
\end{eqnarray}

For any Laurent polynomial $P(z,a)$ of variables $z$ and $a$, we will use $E(P(z,a))$ and $e(P(z,a))$ to denote the highest and lowest power of $a$ in $P(z,a)$. Furthermore, if we write $P(z,a)$ as a polynomial of $a$ with polynomials of $z$ as its coefficients, then the highest power term in the coefficient polynomials of the $a^{E(P(z,a))}$ and the $a^{e(P(z,a))}$ terms are denoted by $p_0^h(P(a,z))$ and $p_0^\ell(P(a,z))$ respectively (these are monomials in the variable $z$). In the case that $P(z,a)=H(D,z,a)$, we will abbreviate $E(P(z,a))$,  $e(P(z,a))$, $p_0^h(P(a,z))$ and $p_0^\ell(P(a,z))$ by $E(D)$, $e(D)$, $p_0^h(D)$ and $p_0^\ell(D)$ respectively. 

\medskip
\begin{example}\label{Example1}{\em
For example, if $D$ is the positive  $T(2,4)$ torus link with its two components assigned anti-parallel orientations, we have $H(D)=-z^{-1}a^{-5}+(z^{-1}+z)a^{-3}+za^{-1}$. Hence $E(D)=-1$, $e(D)=-5$, $p_0^h(D)=z$ and $p_0^\ell(D)=-z^{-1}$. 
}
\end{example}

\medskip
Given an oriented link diagram $D$, let $w(D)$ be the writhe of $D$ and $s(D)$ be the number of Seifert circles in $S(D)$. The following result is well known.

\begin{proposition}\cite{FW,LDH,Morton1986, Ya}
Let $D$ be any link diagram, then $E(D)\le s(D)-w(D)-1$ and $e(D)\ge -s(D)-w(D)+1$. It follows that $s(D)\ge (E(D)-e(D))/2+1$, hence $\textbf{b}(D)\ge \xi(D)= (E(D)-e(D))/2+1$ where $\textbf{b}$ is the braid index of $D$.
\label{thm:MWF}
\end{proposition}

\begin{corollary}\label{Cor4.1}
Let $D$ be any link diagram, then $E(D)\le s(D)-w(D)-1-2r^-(D)$ and $e(D)\ge -s(D)-w(D)+1+2r^+(D)$. 
\end{corollary}

\begin{proof} 
By Remark \ref{reduction_number}, $D$ is equivalent to $D^\p$ with $s(D^\p)=s(D)-r^+(D)$ and $w(D^\p)=w(D)-r^+(D)$, thus we have $e(D)=e(D^\p)\ge -s(D^\p)-w(D^\p)+1=-s(D)-w(D)+1+2r^+(D)$. Similarly,
$E(D)=E(D^\pp)\le s(D^\pp)-w(D^\pp)-1=s(D)-w(D)-1-2r^-(D)$. 
\end{proof}

In the rest of this paper, we will use $c(D)$ to denote the number of crossings in the diagram $D$. If $D$ is a positive diagram (meaning the crossings in $D$ are all positive), then $c(D)=w(D)$. The following result is known, it can be proved using the two special link trivialization operations called Operations P and N in \cite{DHL2020}. Note that it is applicable to any link diagram whose crossings are all positive, the link diagrams do not have to be alternating.

\begin{proposition}\label{E_positive}
Let $D$ be a positive link diagram (that is, all crossings in $D$ are positive), then $E(D)=s(D)-c(D)-1$ and $p_0^h(D)=z^{c(D)-s(D)+1}$.
\end{proposition}

\medskip
We shall also need the following two well known equalities regarding the HOMFLY-PT polynomial. 

\begin{proposition}\label{connect_prop}\cite{HOMFLY}
Let $K_1$ and $K_2$ be two links and let $K_1\#K_2$, $K_1\sqcup K_2$ be the connected sum and disjoint sum of $K_1$ and $K_2$ respectively, then 
$$
H(K_1\#K_2)=H(K_1)H(K_2)
$$ 
and 
$$
H(K_1\sqcup K_2)=H(K_1)H(K_2)\left(\frac{a-a^{-1}}{z}\right).
$$
\end{proposition}

\section{The proof of the main result}\label{proof_sec}

We divide the proof of the theorem into several subsections to make it easier for our reader to follow. Furthermore, we will only consider the case when the $D$ is a positive diagram. If $D$ is a negative diagram, we can consider its mirror image instead since the braid index does not change when we pass from a link to its mirror image. Our goal is to derive explicit formulas for $E(\D_k)$ and $e(\D_k)$, from which we can then calculate $\xi(\D_k)=(E(\D_k)-e(\D_k))/2+1$ and apply Proposition \ref{thm:MWF}.

\subsection{The Seifert graph of  $\D$}\label{sub_reconstruction}

\begin{figure}[h!]
\includegraphics[scale=.8]{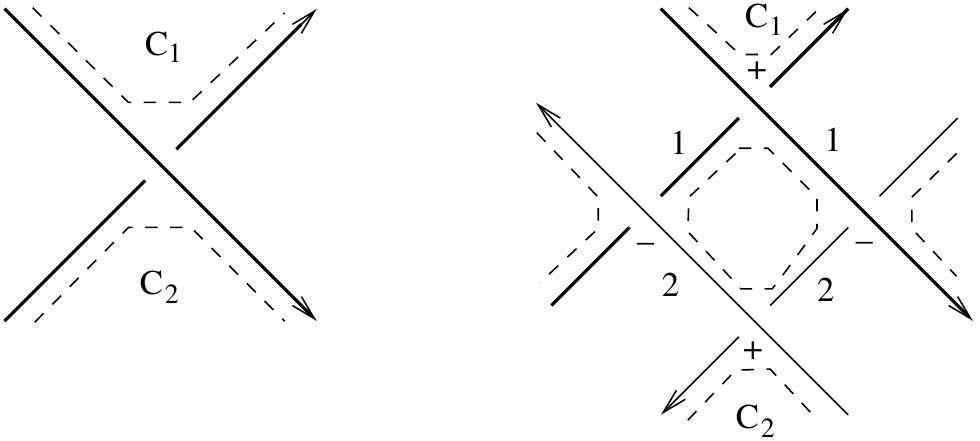}
\caption{Left: A crossing of $D$ between two Seifert circles $C_{1}$ and $C_{2}$ as marked by dashed lines; Right: The four crossings of $D$ and its parallel copy (marked by thin lines) corresponding to the crossing at left. The $s$-circle in  the  middle is a medium $s$-circle while the ones to its right and left are small $s$-circles.  }
\label{Dconstr}
\end{figure}

We shall take a systematic approach to obtain $G_S(\D)$ and to gain a good understanding  of it. Let us choose a parallel copy $D^{\p}$ of $D$ so that its strand is always on our right hand side as we travel on the strand of $D$ following its orientation. Figure \ref{Dconstr} shows the details at a crossing of $D$ and the corresponding crossing junction in $\D$. The crossing is positive and must be between two distinct Seifert circles $C_{1}$ and $C_{2}$ as marked by dashed lines in the figure. Since $D$ and $D^{\p}$  have different orientations, the two positive crossings are between strands from the same component while the negative crossings involve strands from  different components. After smoothing the crossings, it  is obvious that $C_{1}$ and $C_{2}$ remain $s$-circles, with another $s$-circle inserted between them and this new $s$-circle shares a lone positive crossing both  with $C_{1}$ and $C_{2}$. Extending this to all other  4 crossing-junctions of $\D$, we see that in general all $s$-circles with anti-clockwise orientation consist of strands of $D$ while all $s$-circles with clockwise orientation consist of strands of $D^\p$.
Furthermore, $s$-circles of $S(\D)$ can be divided into three groups. The first group contains the original $s$-circles from $D$ which we shall call {\em large $s$-circles}. The second group contains the $s$-circles created in the middle of each 4 crossing-junction of $\D$ as shown in Figure \ref{Dconstr}. Each such $s$-circle  corresponds to a crossing of $D$ and shares lone crossings with two and only two large $s$-circles. We shall call these $s$-circles {\em medium $s$-circles}). The third group contains the remaining $s$-circles, which are obtained   by deleting the negative crossings in $\D$ and  we shall call these {\em small $s$-circles}. Let us  call the vertices in $G_{S}(\D)$ corresponding to  large, medium and small $s$-circles {\em large, medium} and {\em small vertices} respectively.  Notice that each medium $s$-circle also shares  lone negative crossings with exactly two small $s$-circles.

\begin{figure}[htb!]
\includegraphics[scale=.7]{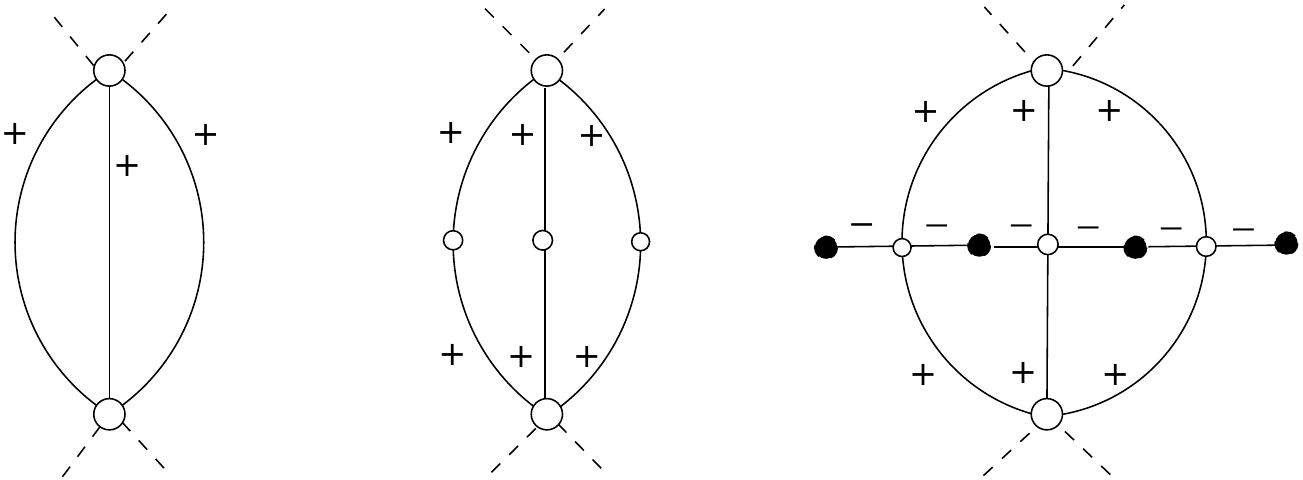}
\caption{Left: Two vertices in $G_S(D)$ connected by three positive edges; Middle: A medium vertex (marked by small circle) is added to the middle of each edge; Right: A small vertex is added to each face and an edge is added between a small vertex and any medium vertex  on the  boundary of this face.}
\label{Figure1}
\end{figure}

\medskip
Let us denote by $\DD$ the (positive) diagram that contains the large $s$-circles  and the medium $s$-circles. $G_{S}(\DD)$ and $G_{S}(\D)$ can be constructed  from $G_S(D)$ in the following way. First we insert a medium vertex to the middle of each edge of $G_S(D)$. This splits each edge of $G_S(D)$ into two lone edges and the resulting graph is $G_{S}(\DD)$. Next we place a small vertex in each  face  of $G_{S}(\DD)$, and add a lone edge (corresponding to a lone negative crossing in $\D$) connecting this vertex  to each medium vertex on the boundary  of this face. The result is $G_{S}(\D)$. See Figure \ref{Figure1}.

\begin{remark}{\em
We make the following observations. First, the length of any cycle in $G_{S}(\DD)$ is a multiple of 4, while the length of any cycle in $\G_{S}(\DD)$ is a multiple of 8. $G_{S}(\D)$  has $s(D)$ large vertices and $c(D)$ medium vertices. The total number of small vertices of $G_{S}(\D)$  is the same as the total number of faces of $G_{S}(D)$ (since a face of $G_{S}(\DD)$ is also a face of $G_{S}(D)$),  which is 
$c(D)-s(D)+2.$ Thus the total number of vertices in $G_{S}(\D)$ (which is the same as $s(\D)$) is  
$
s(D)+c(D)+c(D)-s(D)+2=2c(D)+2.
$
Also, $G_{S}(\D)$ has $4c(D)$ edges and $w(\D)=0$.
}
\end{remark}

\subsection{The determination of $e(\DD)$ and $p_0^\ell(\DD)$.} In  this subsection we  prove the following result.

\medskip
\begin{lemma}\label{DD_lemma} For any reduced alternating knot diagram $D$ whose crossings are all positive, we have
$e(\DD)=s(D)-3c(D)-1$ and $p_0^\ell(\DD)=(-z)^{s(D)-c(D)-1}$.
\end{lemma}

\begin{proof}
Consider the class $\A$ of positive alternating link diagrams that are formed in the following way. Each diagram in $\A$ contains a set of $s$-circles called {\em large $s$-circles}, and a set of $s$-circles called {\em medium $s$-circles} satisfying the following conditions: (1) These $s$-circles do not contain each other in their interiors; (2) Each medium $s$-circle shares lone (positive) crossings with two and only two large $s$-circles; (3) The diagram is reduced. Let $\tilde{\DD}$ be such a link diagram with $s_1(\tilde{\DD})$ being the number of large $s$-circles in $\tilde{\DD}$ and $s_2(\tilde{\DD})$ being the number of medium $s$-circles in $\tilde{\DD}$ (so $s(\tilde{\DD})=s_1(\tilde{\DD})+s_2(\tilde{\DD})$ and $w(\tilde{\DD})=c(\tilde{\DD})=2s_2(\tilde{\DD})$). We claim that $e(\tilde{\DD})=s_1(\tilde{\DD})-3s_2(\tilde{\DD})-1$ and $p_0^\ell(\tilde{\DD})=(-z)^{s_1(\tilde{\DD})-s_2(\tilde{\DD})-1}$. Since $\DD$ belongs to $\A$ with $s_1(\DD)=s(D)$, $s_2(\DD)=c(D)$, the statement of the lemma is proved if we can prove this claim.

\medskip
Use induction on $s_2(\tilde{\DD})\ge 2$. When $s_2(\tilde{\DD})= 2$, it is necessary that $s_1(\tilde{\DD})=2$ as well, so $\tilde{\DD}$ is the link diagram given in Example \ref{Example1} with $e(\tilde{\DD})=-5$ and $p_0^\ell(\tilde{\DD})=-z^{-1}$, hence the statement of the claim holds.

Assume that the statement is true for $s_1(\tilde{\DD})=n\ge 2$. Consider the case when we have $s_1(\tilde{\DD})=n+1$ medium $s$-circles. Consider a medium $s$-circle $\alpha$ that shares lone crossings $x_1$ and $x_2$ with two large $s$-circles $C_1$ and $C_2$. 

Case 1. $\alpha$ is the only medium $s$-circle between its two neighboring large $s$-circles. In this case $x_1$ (and also $x_2$) corresponds to a single edge in $G_S(\tilde{\DD})$ which cannot be a bridge edge. Apply (\ref{Skein1}) to $x_1$.
Observe that $\tilde{\DD}_-$ simplifies (by a Reidemeister II move) to a new diagram $\tilde{\DD}_-^{\p}$ in $\A$ with $s_1(\tilde{\DD}_-^{\p})=s_1(\tilde{\DD})-1$, $s_2(\tilde{\DD}_-^{\p})=s_2(\tilde{\DD})-1$. On the other hand, the smoothing of $x_1$ may result in some additional bridge edge pairs in $\tilde{\DD}_0$ other than $x_2$. Let $k$ be the number of such edge pairs which correspond to nugatory pairs of crossings in $\tilde{\DD}_0$. The removing of these nugatory crossings (including $x_{2}$) yields a link diagram $\tilde{\DD}_0^{\p}$ that is sill in $\A$ with $s_1(\tilde{\DD}_0^{\p})=s_1(\tilde{\DD}_{0}^{\p})-k$, $s_2(\tilde{\DD}_0^{\p})=s_2(\tilde{\DD})-1-k$. By the induction hypothesis we now have 
\begin{eqnarray*}
-2+e(\tilde{\DD}_-)&=&-2+(s_1(\tilde{\DD})-1)-3(s_2(\tilde{\DD})-1)-1=s_1(\tilde{\DD})-3s_2(\tilde{\DD})-1
\end{eqnarray*}
and 
\begin{eqnarray*}
-1+e(\tilde{\DD}_0)&=&-1+(s_1(\tilde{\DD})-k)-3(s_2(\tilde{\DD})-1-k)-1=s_1(\tilde{\DD})-3s_2(\tilde{\DD})+1+2k.
\end{eqnarray*}
Thus $e(\tilde{\DD})=s_1(\tilde{\DD})-3s_2(\tilde{\DD})-1$ since $s_1(\tilde{\DD})-3s_2(\tilde{\DD})+1+2k>s_1(\tilde{\DD})-3s_2(\tilde{\DD})-1$. 

Case 2. $\alpha$ is not the only medium $s$-circle between its two neighboring large $s$-circles. Say there are $m\ge 1$ other medium $s$-circles sharing crossings with $C_{1}$ and $C_{2}$. Again we apply (\ref{Skein1}) to $x_1$. In this case, $\tilde{\DD}_0$ simplifies to $\tilde{\DD}_0^{\p}\in \A$ (after removing the nugatory crossing $x_{2}$) with  $s_1(\tilde{\DD}_0^{\p})=s_1(\tilde{\DD}_{0}^{\p})$, $s_2(\tilde{\DD}_0^{\p})=s_2(\tilde{\DD})-1$. On the other hand, $\tilde{\DD}_-$ simplifies to a diagram that is the connected sum of a diagram $\tilde{\DD}_-^{\p}\in\A$ with $m$ copies of positive $T(2,2)$, and $s_1(\tilde{\DD}_-^{\p})=s_1(\tilde{\DD})-1$, $s_2(\tilde{\DD}_-)=s_2(\tilde{\DD})-1-m$. We have $H(T(2,2))=-z^{-1}a^{-3}+(z^{-1}+z)a^{-1}$. Thus by the induction hypothesis and Proposition \ref{connect_prop} we have
\begin{eqnarray*}
&&-2+e(\tilde{\DD}_-)=-2-3m+e(\tilde{\DD}_-^{\p})\\
&=&-2-3m+(s_1(\tilde{\DD})-1)-3(s_2(\tilde{\DD})-1-m)-1=s_1(\tilde{\DD})-3s_2(\tilde{\DD})-1\\
&<&
-1+e(\tilde{\DD}_0)\\
&=&-1+s_1(\tilde{\DD})-3(s_2(\tilde{\DD})-1)-1=s_1(\tilde{\DD})-3s_2(\tilde{\DD})+1.
\end{eqnarray*}
Thus $e(\tilde{\DD})=s_1(\tilde{\DD})-3s_2(\tilde{\DD})-1$. 

\medskip
Notice that in both cases $e(\tilde{\DD})$ comes from the term $a^{2}H(\tilde{\DD}_-)$. Thus in the first case we have
\begin{eqnarray*}
p_{0}^{\ell}(\tilde{\DD})&=&p_{0}^{\ell}(\tilde{\DD}_-)\\
&=&
(-z)^{(s_{1}(\tilde{\DD})-1)-(s_{2}(\tilde{\DD})-1)-1}\\
&=&(-z)^{s_{1}(\tilde{\DD})-s_{2}(\tilde{\DD})-1},
\end{eqnarray*}
 and in the second case we also have
\begin{eqnarray*}
p_{0}^{\ell}(\tilde{\DD})&=&(-z)^{-m}p_{0}^{\ell}(\tilde{\DD}_-^{\p})\\
&=&
(-z)^{-m+(s_{1}(\tilde{\DD})-1)-(s_{2}(\tilde{\DD})-1-m)-1}\\
&=&(-z)^{s_{1}(\tilde{\DD})-s_{2}(\tilde{\DD})-1}.
\end{eqnarray*}
So the statement of the claim also holds and the lemma is proved.
\end{proof}

\medskip
\subsection{The determination of $E(\D)$ and $p_{0}^{h}(\D)$.} In this subsection we  prove the following lemma.

\medskip
\begin{lemma}\label{E(D)lemma} For any reduced alternating knot diagram $D$ whose crossings are all positive, we have $E(\D)=2s(D)-1$ and $
p_{0}^{h}(\D)=z^{2c(D)-2s(D)+1}.
$
\end{lemma}

Let us call an edge in $G_{S}(\D)$ a {\em positive (negative) edge} if the crossing in $\D$ corresponding to it is positive (negative). Recall that positive edges are the ones connecting a large vertex and a medium vertex while all the negative edges are the ones connecting a medium vertex and a small  vertex.  Each medium vertex has two  negative edges connected to it  and  these two edges connect  to two distinct small vertices so  they are  uniquely associated with this medium vertex. This naturally divides the negative edges of $G_{S}(\D)$, hence the negative crossings of $\D$, into $c(D)$ pairs. As illustrated in Figure \ref{Dconstr}, it is necessary that the two strands at the negative crossings belong to different components. Mark the two components of $\D$ by 1 and 2. For each pair of these negative crossings, the strand belonging to component 1 is the under strand at one and only one of the two crossings (as shown in Figure \ref{Dconstr}), and we will choose this crossing and place it in a crossing set $\C$. That is, $\C$ contains $c(D)$ negative crossings and at  these crossings, the strands belonging to component 1 are always the under strands. Let $\C^\p$ be the set of negative edges in $G_{S}(\D)$ corresponding to the crossings in $\C$. 

\begin{figure}[htb!]
\includegraphics[scale=.9]{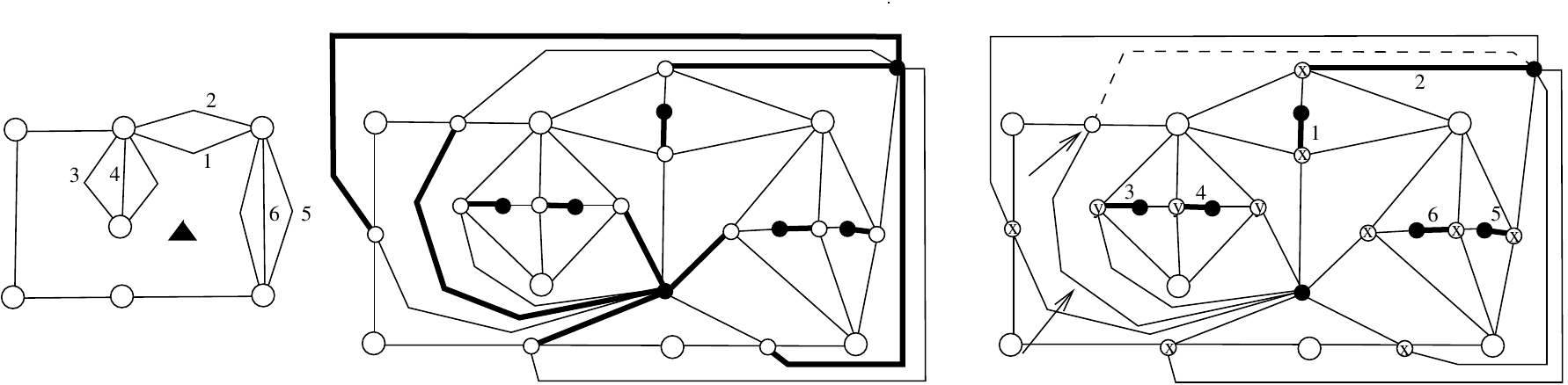}
\caption{Left: $G_{S}(D)$ of a special diagram $D$ where the edges numbered are the proper edges by the order they are deleted in obtaining the spanning tree of $G_{S}(D)$. The initial face is marked by a triangle; Middle: How the negative edges in $\C^\p$ are paired with the negative edges not in $\C^\p$ (marked by thick lines). The medium $s$-circles are marked by small circles and the small $s$-circles are marked by solid dots; Right: The negative edges not in $\C^{\p}$ that are associated with the proper edges of $G_{S}(D)$ (as indicated in the figure at left)  are highlighted by thick lines. }
\label{choice}
\end{figure}

\medskip
Furthermore, if we go  around a small vertex in $G_{S}(\D)$ in any orientation, we encounter the edges from $\C^\p$ and not from $\C^\p$ alternately, see Figure \ref{choice} for an example. This can be observed by following the strand of a small $s$-circle, since the crossings we encounter along this $s$-circle are always negative, we either always arrive at the crossings on the overpass or always arrive at the crossings on the underpass (depending on the orientation we choose to travel), but the strands have to alternate between component 1 and component 2. 

\medskip 
Consider an edge $\gamma$ of $G_S(D)$ that is on the boundary $\partial F$ of a face $F$ in $G_S(D)$. If the lone edge connecting the small vertex of $G_{S}(\D)$ that is placed in $F$ and the medium vertex of $G_{S}(\D)$ that is on $\gamma$ belongs to $\C^\p$, we will say that $\gamma$ is {\em proper} with respect to $F$, otherwise we say that $\gamma$ is {\em improper} with respect to $F$. As we travers $\partial F$ in any given orientation, we encounter proper and improper edges of $G_S(D)$ with respect to $F$ alternately. We now construct a spanning tree of $G_{S}(D)$ in the following way.
 
 \medskip
 Step 1. Start with any face $F_{0}$ of $G_{S}(D)$, and choose any proper edge $\gamma_{1}$ with respect to $F_{0}$, and delete $\gamma_{1}$ from the graph. This eliminates the face $F_{0}$.
 
 \medskip
 Step 2. Let $F_{1}$ be the face that shares $\gamma_{1}$ with $F_{0}$ on its boundary, and choose a proper edge $\gamma_{2}$ with respect to $F_{1}$. Notice that it  is necessary that $\gamma_{2}\not=\gamma_{1}$. Delete  $\gamma_{2}$ from the graph and this eliminates the face $F_{1}$. 
 
 \medskip
Steps 3 to $\beta=c(D)-s(D)+1$. We now continue this process. At each step, we choose a proper edge from the current face, which shares the proper edge chosen to be deleted from the previous edge on its boundary, to be deleted from the graph. Since there are $c(D)-s(D)+2$ faces, we can do this $\beta=\beta(D)=c(D)-s(D)+1$ times and at the end we reach a spanning tree of $G_{S}(D)$. Notice that each $\gamma_{j}$ corresponds to a unique negative edge  $\gamma_{j}^{\p}$ not belonging to $\C^{\p}$, see the right of Figure \ref{choice}. Notice that on the right side of Figure \ref{choice}, vertices marked by the same letter ($x$ or $y$) are contracted to the same vertex due to the MP operations which are always performed in the direction from a medium vertex to a small vertex. The medium vertex pointed by the arrow is also contracted to the same vertex marked by $x$.

\medskip
\begin{remark}\label{betaremark}{\em
We observe first that $\beta=c(D)-s(D)+1$ MP operations can be performed on the crossings corresponding to the $\gamma_{j}^{\p}$'s, in the sequential order of $\gamma_{\beta}^{\p}$, $\gamma_{\beta-1}^{\p}$, ..., $\gamma_{2}^{\p}$, $\gamma_{1}^{\p}$.  The reason is that at each step $\gamma_{j}^{\p}$ brings with it a new medium vertex corresponding to a medium $s$-circle that has not been affected by the previous MP operations, hence $\gamma_{j}^{\p}$ corresponds to a negative lone crossing and an MP operation can be performed on it using the direction from the medium $s$-circle to the small $s$-circle. Now consider the case when an edge $\gamma_{0}\in \C^{\p}$ is deleted (meaning its corresponding negative crossing is smoothed). Let $\gamma_{0}^{\p}\not\in \C^{\p}$ be the edge that is paired with $\gamma_{0}$ and let $F_{0}$ be the face of $G_{S}(\D)$ that contains $\gamma_{0}^{\p}$, and use $F_{0}$ as the starting face to obtain the negative edges $\gamma_{j}^{\p}$, $1\le j\le \beta$. This time, as we perform the MP operations, the medium vertex that $\gamma_{0}^{\p}$ is connected to is shielded away from being swallowed by the combined vertices (due to the MP operations) since $\gamma_{0}$ has been deleted. See the right side of Figure \ref{choice} for an illustration of this, where the edge marked by dashed line is deleted hence the vertex pointed by the arrow is not contracted by the $\beta$ MP operations. Therefore $\gamma_{0}^{\p}$ remains a lone edge after the $\beta$ MP operations discussed above are performed, and we can perform one more MP operation on the crossing corresponds to it. This means that if we smooth a negative crossing in $\C$, then we can perform MP operations on $\beta+1=c(D)-s(D)+2$ negative crossings, regardless whether some crossings in $\C$ have been flipped or not, since the MP operations in the above discussion only used negative crossings not in $\C$.
 }
 \end{remark}
 
 We are now ready to prove Lemma \ref{E(D)lemma}.

 \begin{proof} Choose any crossing in $\C$ and apply skein relation (\ref{Skein2}) to it. We have 
 $H(\D)=H(D_{-})=a^{2}H(D_{+})-azH(D_{0})$. Notice that $s(D_{0})=s(\D)=2c(D)+2$, $w(D_{0})=w(\D)+1=1$. By Corollary \ref{Cor4.1} and Remark \ref{betaremark}, we have
 \begin{eqnarray*}
 E\left(-azH(D_{0})\right)
 &\le&
  1+s(D_{0})-w(D_{0})-1-2r^{-}(D_{0})\\
  &=& 1+2c(D)+2-1-1-2r^{-}(D_{0})\\
 & \le& 
 1+2c(D)-2(c(D)-s(D)+2)\\
 &=&
 2s(D)-3.
 \end{eqnarray*}
 Thus the $-azH(D_{0})$ term will not make a contribution to the $a^{2s(D)-1}$ term. We now consider $D_{+}$, which is  obtained by flipping the crossing. We will choose another crossing from $\C$ and repeat the above argument. Each time we can ignore the one obtained by smoothing the crossing. At the end, we arrive at a single term of the form $a^{2c(D)}H(D^{f})$ where the diagram $D^{f}$ is obtained by flipping all crossings in $\C$.
We observe that $D^{f}$ separates into two disjoint copies of $D$ (since component 1 now sits on top of component 2 at all crossings where they intersect), hence $H(D^{f})=H^{2}(D)(a-a^{-1})/z$ by Proposition \ref{connect_prop}. By Proposition \ref{E_positive}, we have (keep in mind that $w(D)=c(D)$):
$$
E\left(a^{2c(D)}H(D^{f})\right)=2c(D)+2s(D)-2w(D)-2+1=2s(D)-1,
$$ 
and 
$$
p_{0}^{h}\left(a^{2c(D)}H(D^{f})\right)=z^{2c(D)-2s(D)+2-1}=z^{2c(D)-2s(D)+1}.
$$
This proves that $E(\D)=2s(D)-1$ and $
p_{0}^{h}(\D)=z^{2c(D)-2s(D)+1}.
$
\end{proof}

\medskip
\subsection{The determination of $e(\D)$ and $p_{0}^{\ell}(\D)$.} In this section we prove the following lemma.

\begin{lemma}\label{e(D)lemma}
If $D$ is a special alternating knot diagram with positive crossings, then
$e(\D)=2s(D)-2c(D)-3$
and 
$
p_{0}^{\ell}(\D)=-z^{2s(D)-3}.
$
\end{lemma}

\begin{definition}\label{genDdef}{\em
Let $\tilde{D}\in \A$, where $\A$ is the class of alternating knot diagrams defined in the proof of Lemma \ref{DD_lemma}, such that the graph $G^{\p}$ obtained from $G_{S}(\tilde{D})$ by removing the medium vertices is proper. Consider the link diagram corresponding to a Seifert graph  obtained by placing a vertex (called a {\em small vertex}) in each face of $G_{S}(\tilde{D})$, then adding an arbitrary number of edges between such a small vertex and the medium vertices on the boundary of the face that contains the small vertex. Each such edge corresponds to a negative crossing in the corresponding diagram. If a total of $k\ge 0$ such edges are added, then the corresponding link diagram is denoted by $\hat{D}_{k}$.
}
\end{definition}

\begin{lemma}\label{general_e(D)lemma}
Let $\hat{D}_{k}$ be as defined in Definition \ref{genDdef}. Let $s_{1}$ and $s_{2}$ be the numbers of large and medium vertices in $G_{S}(\tilde{D})$, and $f=s_{2}-s_{1}+2$ be the number of faces in $G_{S}(\tilde{D})$ (which is the same as the number of small vertices in $G_{S}(\hat{D}_{k})$), then 
$e(\hat{D})=k+2s_{1}-4s_{2}-3$
and 
$
p_{0}^{\ell}(\hat{D})=(-z)^{k+2s_{1}-2s_{2}-3}.
$
\end{lemma}

\begin{figure}[htb!]
\includegraphics[scale=0.9]{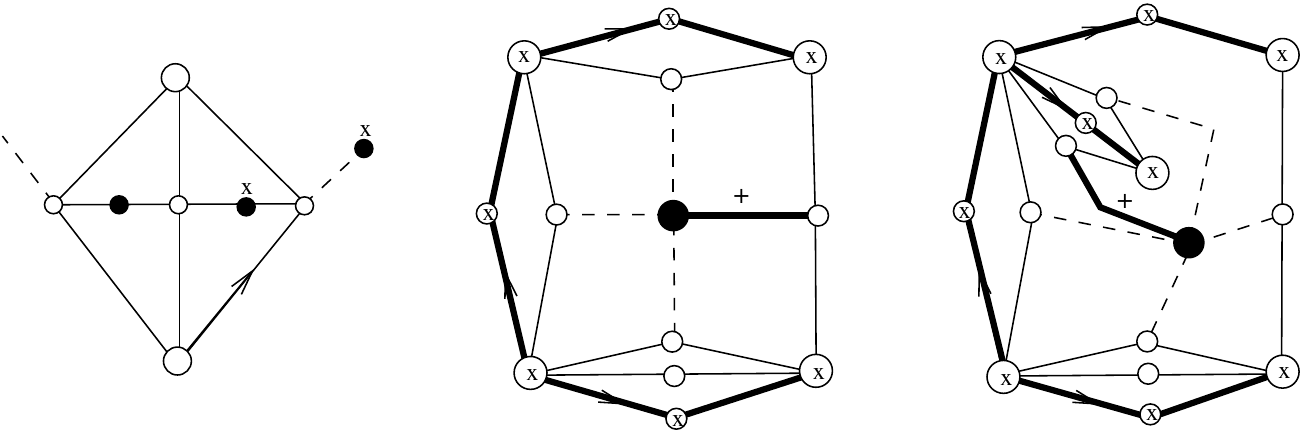}
\caption{Left: The effect of performing an MP operation on a positive lone crossing in $\hat{D}_{k}$. The edge used here corresponds to an edge in the spanning tree $T_{F}$;  Middle and Right:  Two examples of how  the edges of the spanning tree $T_{F}$ are used to perform the MP operations without affecting the boundary of the face that contains $v_{s}$ and the edge corresponding to the flipped crossing. The arrows indicate the directions of the corresponding MP operations and the vertices marked by the same letter are contracted to the same vertex due to the MP operations. The marked positive edge remains a lone edge.}
\label{e(D)fig2}
\end{figure}

\begin{proof}
We use induction on $k$. If $k=0$, then $\hat{D}$ is a disjoint union of $\tilde{D}$ and $f$ copies of trivial knots, hence the result follows from Proposition \ref{connect_prop} and Lemma \ref{DD_lemma}. Assume that this result is true for $k\ge 0$ and consider the case $k+1$. Choose any negative crossing in $\hat{D}_{k+1}$ and apply (\ref{Skein2}) to it: $H(D_{-})=a^{2}H(D_{+})-azH(D_{0})$ with $D_{-}=\hat{D}_{k+1}$ and $D_{0}=\hat{D}_{k}$ (to which the induction hypothesis applies). There are two cases to consider. 

Case 1: The negative crossing is not a lone crossing. In this case it is necessary that $k+1\ge 2$ and $H(D_{+})$ simplifies (via a Reidemeister II move) to $\hat{D}_{k-1}$. Thus we have 
$$
e(a^{2}H(D_{+}))=2+e(\hat{D}_{k-1})=2+k-1+2s_{1}-4s_{2}-3=(k+1)+2s_{1}-4s_{2}-3,
$$
with 
$$
p_{0}^{\ell}(a^{2}H(D_{+}))=p_{0}^{\ell}(\hat{D}_{k-1})=(-z)^{(k+1)+2s_{1}-2s_{2}-5}.
$$
On the other hand, we have
$$
e(-azH(D_{0}))=1+e(\hat{D}_{k})=1+k+2s_{1}-4s_{2}-3=(k+1)+2s_{1}-4s_{2}-3,
$$
with 
$$
p_{0}^{\ell}((-az)H(D_{0}))=(-az)p_{0}^{\ell}(\hat{D}_{k})=(-z)^{(k+1)+2s_{1}-2s_{2}-3}.
$$
It follows that $e(\hat{D}_{k+1})=(k+1)+2s_{1}-4s_{2}-3$ and $p_{0}^{\ell}(\hat{D}_{k})=(-z)^{(k+1)+2s_{1}-2s_{2}-3}$.

Case 2. The negative crossing is a lone crossing, which corresponds to a (negative) lone edge connecting a small vertex $v_{s}$ and a medium vertex $v_{m}$ in $G_{S}(\hat{D}_{k+1})$. Again in this case the $-azH(D_{0})$ term yields the needed $e(\hat{D}_{k+1})$ and $p_{0}^{\ell}(\hat{D}_{k})$, so it suffices  to show that $e(a^{2}H(D_{+}))>(k+1)+2s_{1}-4s_{2}-3$. Let $F$ be the face of $G_{S}(\hat{D}_{k+1})$ containing $v_{s}$ and let $T_{F}$ be a spanning tree of $G^{\p}$ that does not contain any edge on the boundary of $F$. Since any edge in $T_{F}$ corresponds to two large vertices each connected to a common medium vertex by a (positive) lone edge, an MP operation can be performed using one of these lone crossings as shown in the left of Figure \ref{e(D)fig2}. The effect of this move is that $C_{1}$ and $C_{2}$ are swallowed by the newly created $s$-circle. Since there are $s_{1}-1$ edges in $T_{F}$, we can perform $s_{1}-1$ such operations. Furthermore, these operations do not affect the medium vertices on the boundary of $F$, hence the flipped crossing remains a positive lone crossing, and one more MP operation can be performed. Thus we have $r^{+}(D_{+})\ge s_{1}$. Notice that $s(D_{+})=2s_{2}+2$ and $w(D_{+})=2s_{2}+2-(k+1)$. It follows from Corollary \ref{Cor4.1} that
\begin{eqnarray*}
e(a^{2}H(D_{+}))&\ge& 2-s(D_{+})-w(D_{+})+1+2r^{+}(D_{+})\\
&\ge &2-2(s_{2}+2)-(2s_{2}+1-k)+1+2s_{1}\\
&=&k+2s_{1}-4s_{2}>(k+1)+2s_{1}-4s_{2}-3.
\end{eqnarray*}
 \end{proof}
  
 We are now ready to prove Lemma \ref{e(D)lemma}.
 
\begin{proof}
If all maximal paths in $\G_{S}(D)$ are proper, then the statement of the lemma follows from Lemma \ref{general_e(D)lemma} by substituting $s_{1}=s(D)$, $s_{2}=c(D)$ and $k=2c(D)$. Let us consider the case that $\G_{S}(D)$ contains maximal paths that are not proper. Remark \ref{maximalpathremark}, these maximal paths contain vertices that have degree 2 in $G_{S}(D)$. Furthermore, if there are multiple vertices with this property in a maximal path, they appear consecutively along the path, and the two end vertices of the path musth have degree more than  2 in  $G_{S}(D)$. Let $k_0$ be the total number of degree 2 vertices in $G_{S}(D)$. Each such vertex is a large vertex of $G_{S}(\DD)$  which has two lone edges connected to it. Of these two edges, we shall choose one that we encounter first by walking through the maximal path (it do  not matter in which direction). Denote the set  of these chosen edges by $X$. We  have $|X|=k_0$. We  now choose any crossing in  $X$  and  apply  (\ref{Skein1}) to it:
$H(D_{+})=a^{-2}H(D_{-})+a^{-1}zH(D_{0})$. We observe  that smoothing a crossing in $X$ increases the number of  MP operations in the diagram by one since the we have created a positive nugatory crossing which corresponds to an  isolated medium $s$-circle. That is, $s(D_{0})=s(\D)=2c(D)+2$, $w(D_{0})=w(\D)-1=-1$, and $r^{+}(D_{0})\ge s(D)$. It follows that 
$$
e(a^{-1}zH(D_{0}))\ge -1 -(2c(D)+2)-(-1)+1+2s(D)=2s(D)-2c(D)-1>2s(D)-2c(D)-3.
$$
Thus we shall ignore the $a^{-1}zH(D_{0})$ term, and will continue the process with the $a^{-2}H(D_{-})$ term. That is, we start with diagram obtained by flipping the first crossing, and choose another crossing in $X$, and apply (\ref{Skein1}) to it. Flipping the previously chosen crossings do not affect our ability to perform the $s(D)-1$ MP operations based on any spanning tree of $G_{S}(D)$ as described in Case 2 of the proof of Lemma \ref{general_e(D)lemma}, thus the same argument always applies to the diagram obtained by smoothing the crossing. At the end, we arrive at a term of the form $a^{-2k_{0}}H(D^{X})$, where $D^{X}$ is obtained by flipping every crossing in $X$. Now observe that for each flipped crossing, a Reidemeister II move can be used to simplify the diagram: the result is that a large $s$-circle and two medium $s$-circles are combined into one medium $s$-circle and the edges connecting to two medium $s$-circles to small vertices are all connected to  this newly created medium vertex. Thus, $D^{X}$ simplifies to a diagram $\hat{D}_{k}$ as defined in Definition \ref{genDdef} with $s_{1}=s(D)-k_{0}$, $s_{2}=c(D)-k_{0}$ and $k=2c(D)$. By Lemma \ref{general_e(D)lemma} we have
\begin{eqnarray*}
&&e\left(a^{-2k_{0}}H(D^{X}))\right)\\
&=& -2k_{0}+e(\hat{D}_{k})\\
&=&
-2k_{0}+2c(D)+2(s(D)-k_{0})-4(c(D)-k_{0})-3\\
&=&
2s(D)-2c(D)-3,
\end{eqnarray*}
and
\begin{eqnarray*}
&&p_{0}^{\ell}\left(a^{-2k_{0}}H(D^{X}))\right)\\
&=& p_{0}^{\ell}\left(\hat{D}_{k}\right)\\
&=&
(-z)^{2c(D)+2(s(D)-k_{0})-2(c(D)-k_{0})-3}\\
&=&
-z^{2s(D)-3}.
\end{eqnarray*}
This shows that $e(\D)=2s(D)-2c(D)-3$ and $p_{0}^{\ell}(\D)=-z^{2s(D)-3}$.
\end{proof}

\subsection{The determination of $E(\D_{k})$ and $e(\D_{k})$.} We are now ready to derive explicit formulas for $E(\D_{k})$ and $e(\D_{k})$ as stated in the following lemma.

\medskip
\begin{lemma}\label{Dklemma}
If $k\ge 0$, then
$e(\D_{k})=2s(D)-2c(D)-2k-3$ and 
$$
E(\D_{k})=\left\{
\begin{array}{ll}
2s(D)-2k-1& \ {\rm if}\ 0\le k\le s(D);\\
-1 &  \ {\rm if}\ k> s(D)
\end{array}
\right.
$$
with 
$$
p_{0}^{h}(\D_{k})=\left\{
\begin{array}{ll}
z^{2c(D)-2s(D)+1}& \ {\rm if}\ 0\le k\le s(D);\\
z &  \ {\rm if}\ k> s(D).
\end{array}
\right.
$$
On the other hand, if $k\le 0$, then
$E(\D_{k})=2s(D)-2k-1$, and 
$$
e(\D_{k})=\left\{
\begin{array}{ll}
2s(D)-2k-2c(D)-3& \ {\rm if}\ k\ge s(D)-c(D)-2;\\
1 &  \ {\rm if}\ k< s(D)-c(D)-2
\end{array}
\right.
$$
with 
$$
p_{0}^{\ell}(\D_{k})=\left\{
\begin{array}{ll}
-z^{2s(D)-3}& \ {\rm if}\ k\ge s(D)-c(D)-1,\ {\rm or}\ k= s(D)- c(D)-2,\ s(D)>2;\\
-2z &  \ {\rm if}\ k= s(D)-c(D)-2,\ s(D)=2;\\
-z &  \ {\rm if}\ k<s(D)- c(D)-2.
\end{array}
\right.
$$
\end{lemma}

\begin{proof}
First consider the case $k\ge 0$. Use induction on $k$. The case of $k=0$ has already been proved. Assume the lemma is true for some $k\ge 0$, consider the case $k+1$. For $\D_{k+1}$, apply (\ref{Skein1}) to one of the positive crossings in the $k+1$ full positive twists: $H(\D_{k})=H(D_{+})=a^{-2}H(D_{-})+a^{-1}zH(D_{0})$. Notice that $D_{0}$ deforms to the trivial knot, while $D_{-}$ simplifies to $\D_{k}$ to which the induction hypothesis applies. Thus we have 
$e(a^{-2}H(D_{-}))=-2+ 2s(D)-2c(D)-2k-3=2s(D)-2c(D)-2(k+1)-3<e(a^{-1}zH(D_{0}))=-1$ for any $k\ge 0$.  On the other hand, if $k+1< s(D)$, then $E(a^{-2}H(D_{-}))=-2+2s(D)-2k-1=2s(D)-2(k+1)-1>E(a^{-1}zH(D_{0}))=-1$, with $p_{0}^{h}(a^{-2}H(D_{-}))=p_{0}^{h}(D_{-})=z^{2c(D)-2s(D)+1}$. If $k+1=s(D)$, then $E(a^{-2}H(D_{-}))=2s(D)-2(k+1)-1=-1=E(a^{-1}zH(D_{0}))$. We have $p_{0}^{h}(a^{-1}zH(D_{0}))=z$ and $p_{0}^{h}(a^{-2}H(D_{-}))=z^{2c(D)-2s(D)+1}$.
Since $c(D)>s(D)$, $2c(D)-2s(D)+1>1$ so $p_{0}^{h}(\D_{k+1})=z^{2c(D)-2s(D)+1}$. If $k+1\ge s(D)+1$, then we have $E(a^{-2}H(D_{-}))=-2-1<E(a^{-1}zH(D_{0}))=-1$
hence $E(\D_{k+1})=-1$ with $p_{0}^{h}(\D_{k+1})=p_{0}^{h}(a^{-1}zH(D_{0}))=z$. This proves the case for $k\ge 0$. 

For $\D_{k}$ with $k<0$, apply (\ref{Skein2}) to one of the negative crossings in the $k$ full negative twists: $H(\D_{k})=H(D_{-})=a^{2}H(D_{+})-azH(D_{0})$. Again $D_{0}$ deforms to the trivial knot, while for $k\le -1$, $D_{+}$ simplifies to $\D_{k+1}$ to which the induction hypothesis applies. Thus we have 
\begin{eqnarray*}
E(\D_{k})&=&E(a^{2}H(D_{+}))\\
&=&2+ 2s(D)-2(k+1)-1\\
&=&2s(D)-2k-1>E(-azH(D_{0}))=1
\end{eqnarray*} 
for any $k\le 0$.  On the other hand, if $s(D)-c(D)-1\le k\le -1$, then 
\begin{eqnarray*}
e(\D_{k})&=&e(a^{2}H(D_{+}))\\
&=&2+2s(D)-2(k+1)-2c(D)-3\\
&=&-2k-2(c(D)-s(D)+1)-1<e(-azH(D_{0}))=1,
\end{eqnarray*} 
with 
$$
p_{0}^{\ell}(\D_{-k})=p_{0}^{\ell}(a^{2}H(D_{+}))
=p_{0}^{\ell}(D_{+})=-z^{2s(D)-3}.
$$
If $k=s(D)-c(D)-2$, and $s(D)=2$, then 
\begin{eqnarray*}
e(\D_{k})&=&e(a^{2}H(D_{+}))\\
&=&2+2s(D)-2(k+1)-2c(D)-3\\
&=&2+2s(D)-2(s(D)-c(D)-1)-2c(D)-3=1=e(-azH(D_{0})),
\end{eqnarray*}
 with $p_{0}^{\ell}(\D_{k})=p_{0}^{\ell}(a^{2}H(D_{+}))+p_{0}^{\ell}(-azH(D_{0}))=-z-z=-2z$. If $k=s(D)-c(D)-2$, but $s(D)>2$, then $e(\D_{k})=e(a^{2}H(D_{+}))=1$, with $p_{0}^{\ell}(\D_{k})=p_{0}^{\ell}(a^{2}H(D_{+}))=-z^{2s(D)-3}$ since $2s(D)-3>1$. Finally, if $k<s(D)-c(D)-2$, then $e(a^{2}H(D_{+}))=3>1=e(-azH(D_{0})=e(\D_{k})$, with $p_{0}^{\ell}(\D_{k})=p_{0}^{\ell}(-azH(D_{0}))=-z$.
\end{proof}

\medskip
Theorem \ref{main_theorem} now follows from Lemma \ref{Dklemma} easily. For example, 
if $D$ contains only positive crossings and $k\ge 0$, then we have
\begin{eqnarray*}
\xi(\D_k)&=&(E(\D_k)-e(\D_k))/2+1\\
&=&
\left\{
\begin{array}{ll}
c(D)+2 & \ {\rm if}\ k\le s(D);\\
c(D)+2+k-s(D)  & \ {\rm if}\ k> s(D)
\end{array}
\right.\\
&=&
c(D)+2+\rho_k(D).
\end{eqnarray*}
The case $k<0$ can be similarly obtained. If $D$ is a negative diagram, we can apply the above to the mirror image of $D$  first, then take the mirror image again. Passing to the mirror image does not change the number of crossings and the number of Seifert circles in the diagram, so the only change is the sign of $k$.

\section{Application to the ropelengths of knots}\label{ropelength_sec}

\medskip
Let $\K$ be an alternating link. We will call the largest braid index among all braid indices corresponding to different orientation assignments to the components of $\K$ the {\em maximum braid index} of $\K$ and denote it by
$\B(\K)$. It has been shown recently  that $R(\K)\ge a_{0} \B(\K)$ for some constant $a_0>0$ that is independent of $\K$ (in fact $a_{0}>1/14$). 
Thus for a link $\K$, if $\B(\K)$ is proportional to $Cr(\K)$, then $L(\K)$ would be bounded below by a constant multiple of $Cr(\K)$. For example, for the two component torus link $T(2,2n)$, we have $\B(T(2,2n))=n+1$, hence $R(T(2,2n))=O(Cr(T(2,2n)))$. However, for a link with small maximum braid index, this result would be of little help to us. As a consequence of Theorem \ref{main_theorem}, we can prove the following theorem.

\medskip
\begin{theorem}\label{MainT}
Let $\K$ be a special alternating knot, then $R(\K)\ge b_0 Cr(\K)$ for some constant $b_0>0$ that is independent of $\K$. 
\end{theorem}

Let us use a lattice realization of $\K$ as a way to estimate $R(\K)$. Let $\K$ be a link and $K_c$ a realization of $\K$ on the cubic lattice. The length of $K_c$ is denoted by $\eta(K_c)$ and the minimum of $\eta(K_c)$ over all lattice realization $K_c$ of $\K$ is denoted by $\eta_c(\K)$.  The following result shows that $\B(\K)$ bounds $R(\K)$ from below with $\eta_c(\K)$ serving as a bridge.  

\begin{proposition}\label{P1}\cite{Diao2002,Diao2020}
For any lattice realization $K_c$ of the link $\K$, we have $\B(\K)<\eta(K_c)$. Passing to the length minimizer, we have $\B(\K)<\eta_c(\K)<14R(\K)$, that is, $R(\K)>(1/14)\B(\K)$.
\end{proposition}

\medskip
We will now provide the proof of Theorem \ref{MainT} as a consequence of Theorem \ref{main_theorem} and Proposition \ref{P1}. 

\begin{proof}
Let $K_c$ be an arbitrary realization of $\K$ in the cubic lattice and let $\eta(K_c)$ be the length of it. The set $\{\x+t(\frac{1}{2},\frac{1}{2},\frac{1}{2}): \x\in K_c,\ 0\le t\le 1\}$ is an embedding of the annulus $\{1\le x^2+y^2\le 2:\ x, y\in \R\}$ into $\R^3$ with $K_c$ and $K_c^\p=(\frac{1}{2},\frac{1}{2},\frac{1}{2})+K_c$ as its boundary curves. Denote by $\L(K_c,K_c^\p)$ the oriented link formed by $K_c$ and $K_c^\p$ such that they are assigned opposite orientations, then $\L(K_c,K_c^\p)$ is equivalent to $\D_k$ with $k=w(D)+Lk$ where $Lk$ is the linking number of $\L(K_c,K_c^\p)$. By Theorem \ref{main_theorem} we have $Cr(\K)+2\le \b(\D_k)\le \B(\L(K_c,K_c^\p))$. If we double the cubic lattice (so every edge is doubled and a new lattice point is inserted in the middle), then we obtain a realization of 
$\L(K_c,K_c^\p)$ on the cubic lattice whose length is $4\eta(K_c)$. By Proposition \ref{P1}, we have $Cr(\K)+2<4\eta(K_c)$. Since $K_c$ is arbitrary, this leads to $Cr(\K)+2<4\eta_c(\K)<56R(\K)$. That is, $R(\K)>a_0 c(\K)$, where $a_0\ge 1/56$ is a constant which does not depend on $\K$.
\end{proof}

\begin{remark}{\em
Of course, the immediate question one may ask is, can Theorem \ref{main_theorem} be generalized to alternating knots that are not special? We note that the proof of Theorem \ref{main_theorem} relied on the nice structure of the Seifert circle decomposition of $\D$, hence that of $G_{S}(\D)$. It is not clear to the author  how to get around this, but  the result in this paper is indeed very encouraging that the ropelength conjecture may hold in general for all alternating knots and links!
}
\end{remark}

\end{document}